\documentclass{amsproc2000}
\usepackage{latexsym}

\copyrightinfo{2000}{American Mathematical Society}

\newtheorem{thm}{Theorem}[section]

\newtheorem{prop}[thm]{Proposition}

\theoremstyle{definition}
\newtheorem{df}[thm]{Definition}
\newtheorem{ex}[thm]{Example}

\theoremstyle{remark}
\newtheorem{rmk}[thm]{Remark}

\newcommand{\pf}{\noindent{\sc Proof.}\ }

\newcommand{\boom}{\quad\lower3pt\hbox{\vrule height1.1ex width .9ex depth -.2ex}
                    \vskip9pt}

\usepackage{euscript}
\renewcommand{\mathcal}[1]{\EuScript{#1}}

\renewcommand{\phi}{\varphi}
\newcommand{\chigh}{{\raise1.5pt\hbox{$\chi$}}}
\newcommand{\Ga}{\Gamma}

\newcommand{\R}{\mathbb{R}}

\newcommand{\D}{\mathop{\raise0.1ex\hbox{$\mathfrak{D}$}}} 

\newcommand{\End}{\mathop{\rm End}}
\newcommand{\Hom}{\mathop{\rm Hom}}

\newcommand{\gog}{\mathfrak{g}}

\newcommand{\co}{\colon\thinspace}

\newcommand{\isom}{\cong}

\newcommand{\act}{\mathbin{\hbox{$<\kern-.4em\mapstochar\kern.4em$}}}
\newcommand{\ract}{\mathbin{\hbox{$\mapstochar\kern-.3em>$}}}

\renewcommand{\Tilde}{\widetilde}

\newcommand{\Ri}{\overrightarrow}

\newcommand{\surj}{-\!\!\!-\!\!\!-\!\!\!\gg}
\newcommand{\inj}{>\!\!\!-\!\!\!-\!\!\!-\!\!\!>}

\newcommand{\lrah}{\hbox{$\,-\!\!\!-\!\!\!
-\!\!\!-\!\!\!-\!\!\!-\!\!\!-\!\!\!\longrightarrow\,$}}

\begin{document}

\title[DIFFERENTIAL OPERATORS AND ACTIONS OF LIE ALGEBROIDS]{DIFFERENTIAL 
OPERATORS AND \\ ACTIONS OF LIE ALGEBROIDS}

\author{Y. Kosmann-Schwarzbach}
\address{Centre de Math{\'e}matiques, UMR 7640 du CNRS\\
        {\'E}cole Polytechnique\\
        F-91128 \\ Palaiseau, France}
\email{yks@math.polytechnique.fr}
\thanks{}

\author{K. C. H. Mackenzie}
\address{Department of Pure Mathematics\\
        University of Sheffield\\
        Sheffield, S3 7RH\\
        United Kingdom}
\email{K.Mackenzie@sheffield.ac.uk}
\thanks{}

\subjclass[2000]{Primary 17B66. Secondary 22A22, 58D19, 58H05}

\begin{abstract}
We demonstrate that the notions of derivative representation of a Lie
algebra on a vector bundle, of semi-linear representation of a Lie
group on a vector bundle, and related concepts, may be understood in
terms of representations of Lie algebroids and Lie groupoids, and we
indicate how these notions extend to derivative representations of
Lie algebroids and semi-linear representations of Lie groupoids in
general. 
\end{abstract}

\maketitle

\section*{Introduction}

This paper deals with actions on vector bundles.
The first part (Sections \ref{sect:dovb}--\ref{sect:drLa})
deals with the
infinitesimal actions of Lie algebras and Lie algebroids, 
while the second part (Sections \ref{groups} and \ref{sect:gar})
deals with the global actions of Lie groups and Lie groupoids.

When passing from the case of actions on vector spaces to that
of actions on vector bundles, the notions of linear endomorphism
and linear isomorphism admit straightforward generalizations 
which, however, have rarely been spelled out in the literature.
On the global level, it is clear that the analogue of
a linear isomorphism is a vector bundle automorphism, not necessarily
base-preserving. Such an automorphism gives rise
to an isomorphism of the vector space of sections
which has the additional property of \emph{semi-linearity}, {\it i.e.},
when a section is multiplied by a function on the base, the image of
this section is multiplied by the given function composed with the 
inverse of the 
diffeomorphism of the base defined by the vector bundle
automorphism. (This could be viewed as a morphism from the space of
sections to itself, equipped with two module structures.)
Actions of Lie groups on a vector bundle therefore give rise to
representations in its space of sections 
by semi-linear isomorphisms. 

The infinitesimal counterpart of such a vector bundle automorphism can 
be determined by differentiating a one-parameter group of semi-linear
isomorphisms of the space of sections. The result is a \emph{derivative
endomorphism}, a first-order differential operator with the additional
property that its symbol is a scalar multiple of the identity.
Thus the derivative endomorphisms of the space of sections of a vector 
bundle generalize the linear endomorphisms of a vector space, and 
the representations of Lie algebras on the space of sections of a
vector bundle by 
derivative endomorphisms generalize the linear representations.

These basic facts (and their cohomological interpretation, which we do
not repeat here) were described by Kosmann-Schwarzbach in 1976
\cite{Kosmann-Schwarzbach:1976}. Actually, semi-linear
isomorphisms and derivative endomorphisms are both particular cases of what
Jacobson, in 1935, called \emph{pseudo-linear transformations} 
\cite{Jacobson:1935}, \cite{Jacobson:1937}.
But Jacobson's purpose was to extend the theory of elementary
divisors, and the relevance of his work for differential geometry
was, to the best of our knowledge, forgotten.

Meanwhile, the theory of Lie groupoids was created by Ehresmann during the 1950s, 
and their infinitesimal counterpart, the Lie algebroids, were introduced by 
Pradines in 1967 \cite{Pradines:2}. Groupoid theory has been developed by many 
authors, in the context of homotopy theory and $C^*$--algebras, as well as in 
differential geometry, and it has come to differ in substantial respects from 
the elementary theory of groups. However, there is a Lie functor 
which associates to a Lie groupoid, $G$,
a Lie algebroid, denoted $AG$, by a construction which extends in a 
straightforward way the construction of the Lie algebra of a Lie group.  
The resulting Lie theory provides a very broad framework which encompasses 
many of the infinitesimal constructions of differential geometry. 

The concept of action of a Lie groupoid on a fibered manifold goes back to 
Ehresmann \cite{Ehresmann:1957}, who showed that when a groupoid 
acts on a space that is fibered over the base of the groupoid, the 
pullback of the groupoid by the fibration map is, 
in a natural way, a groupoid, called the 
\emph{action groupoid} or \emph{transformation groupoid}. 
The development 
of an abstract theory of Lie algebroids is more recent: in \cite{HigginsM:1990a}
Higgins and Mackenzie provided a basic set of 
algebraic constructions for abstract Lie algebroids modelled on those known in groupoid theory. In particular, they defined the infinitesimal 
action of a Lie algebroid on a fibered manifold and the 
associated \emph{action Lie algebroid} or \emph{transformation Lie algebroid}.
The Lie algebroid of an action Lie groupoid is the action Lie algebroid of the
corresponding infinitesimal action; however the action Lie algebroid of an
infinitesimal Lie algebra action may integrate to a Lie groupoid which is
not necessarily an action groupoid. This may at first seem to be a negative
feature, but in fact the groupoid in such a situation provides a global object
encoding the infinitesimal action even though a global action is not available. 
Here, as is generally the case with abstract Lie algebroids,
determining the conditions under which
a particular construction for Lie algebroids will admit a 
global analogue is a difficult problem. See Dazord \cite{Dazord:1997},
and Moerdijk and Mr\v cun \cite{MoerdijkM} for a systematic 
account of the integration of Lie algebroid actions which includes
many applications.

The aim of this paper is to relate the concept of Lie algebroid action to
the earlier work \cite{Kosmann-Schwarzbach:1976}
on representations by derivative endomorphisms (and, 
correspondingly,
Lie groupoid actions and  
representations by semi-linear isomorphisms.) 
There are several aspects to the relationship. 
First, the derivative endomorphisms of the space of sections of
a vector bundle are the sections of a Lie
algebroid (Theorem \ref{liealgd}). In \cite{Mackenzie:LGLADG}, 
Mackenzie called this Lie algebroid the
\emph{Lie algebroid of covariant differential 
operators}, and denoted it by ${\rm CDO}(E)$, for a vector bundle
$E$. Here, we denote it by ${\mathfrak D}(E)$.
A \emph{derivative representation} of a Lie algebra on a vector bundle
is defined as a morphism 
into the Lie algebra of sections of ${\mathfrak D}(E)$,
while a representation of a Lie algebroid on a vector bundle is
a morphism into this Lie algebroid.
We prove that the derivative representations of Lie algebras 
on vector bundles coincide with the representations of action Lie
algebroids (Theorem \ref{thm:drep}).
Then, with a suitable definition of the 
derivative representations  of a Lie algebroid on a vector bundle,
we show that this one-to-one correspondence can be extended to 
the case of Lie algebroid actions (Theorem \ref{thm:drepoid}).

In Sections \ref{groups} and \ref{sect:gar}, we present global 
versions of the preceding results. We define the 
\emph{semi-linear isomorphisms} of the space of sections of a vector 
bundle. A \emph{semi-linear representation} of a Lie group on a vector 
bundle is defined as a morphism into the group of semi-linear isomorphisms 
of its space of sections, satisfying a smoothness condition, while a 
representation of a Lie groupoid on a vector bundle is a morphism of Lie 
groupoids into the Lie groupoid of all linear isomorphisms from a fiber 
to some (in general different) fiber of the bundle. In Theorem \ref{thm:slrep} 
we prove that the semi-linear representations of Lie groups on vector bundles 
coincide with the representations of their action groupoids. 
For a groupoid acting on a fibered manifold, $F$, there is likewise an
action groupoid; considering an
action of the groupoid on a vector bundle with base $F$, we show in
Proposition \ref{prop:gpdslr} that there is an associated semi-linear
representation which is a group representation of the group of
bisections of the groupoid.

In the Appendix, we recall the definition of the twisted derivations 
of an algebra, and we show that the pseudo-linear endomorphisms of 
the module of sections $\Gamma E$ of a vector bundle $E$ with base $M$
can be defined in terms of twisted derivations of an algebra whose 
underlying vector space is $C^{\infty}(M) \oplus \Gamma E$.
In particular, the derivative endomorphisms can be defined in terms of 
the derivations of this algebra. Similarly, the semi-linear isomorphisms 
can be defined in terms of algebra automorphisms.

Throughout the paper, the manifolds that we consider are assumed to be 
smooth and second countable, and all maps are assumed to be smooth.

\medskip

\noindent{\it Acknowledgments}.
The authors would like to extend their warmest thanks to Mike Prest and
Ted Voronov, the organizers of the workshop on ``Quantization, deformations 
and new homological and categorical methods in mathematical physics''
(Manchester, 2001) which was particularly rich in stimulating
lectures. Although the present paper does not reflect the contents of
the authors' lectures in Manchester, it was the workshop that 
made its realization possible.

\section{Infinitesimal automorphisms of vector bundles}
\label{sect:dovb}

Consider a vector bundle $(E, q, M)$ and let $C^\infty(M)$ be the algebra of 
$\R$--valued functions on $M$. 
We are concerned with first-order differential
operators, $D\co\Ga E\to \Ga E$, 
for which there
exists a vector field $D_M$ on $M$ such that
\begin{equation}                                    \label{eq:D}
D(f\psi) = fD(\psi) + D_M(f)\psi \ ,
\end{equation}
for all $f\in C^\infty(M)$ and all $\psi\in\Ga E$.
Below we explain that operators on modules over a ring
satisfying \eqref{eq:D},
or some variant of it, have been studied under a
variety of names. 
In the case of the module of sections of a vector bundle, 
we shall call such a differential operator 
a \emph{derivative endomorphism} of $\Ga E$.

\begin{df}
\label{df:de}
Let $(E,q,M)$ be a vector bundle. A \emph{derivative endomorphism} 
of $\Ga E$ is an $\R$--linear endomorphism $D$ of 
$\Ga E$ such that there exists an $\R$--linear endomorphism, $D_M$,
of $C^\infty(M)$ satisfying \eqref{eq:D}.
\end{df}

It follows that $D$ is a first-order differential operator. It also
follows that $D_M$ is a derivation of $C^\infty(M)$, and therefore 
corresponds to a vector field on the base manifold $M$. 
The differential operator $D$ is of order $0$ if and only if 
$D_M = 0$.

\begin{ex}\label{example}
\rm Given 
any linear connection $\nabla$ in a vector bundle $E$, and any vector 
field $X$ on the base $M$, the covariant derivation $\nabla_X$ satisfies
(\ref{eq:D}).
\end{ex}

Recall that the \emph{symbol}, $\sigma(D)$,
of a first-order differential operator $D$ is
defined by
$$
\sigma(D)(df) = [D,f] \ , 
$$
where $f $ is any smooth function on $M$ and $[D,f]$ is the commutator of
the operators $D$ and multiplication by $f$.
Condition (\ref{eq:D}) is  
equivalent to the requirement that the symbol 
of $D$, evaluated at any $1$--form 
$\xi$ on $M$ at $x \in M$, be a scalar multiple of the identity of
the fiber $E_x$ at $x$, $\sigma(D)\xi = 
\langle \xi, D_M\rangle \mbox{Id}_{E_x}$. 
Therefore, a first-order differential operator
is a derivative endomorphism of $\Ga E$ if and only if it has
scalar symbol.

Whereas the first-order differential operators on a vector bundle of rank $> 1$
do not constitute a Lie algebra under the commutator, 
a simple computation shows that the commutator 
\begin{equation}
\label{eq:bracket}
[D_1, D_2] = D_1\circ D_2 - D_2\circ D_1
\end{equation}
of two derivative endomorphisms is also a derivative endomorphism, and that
\begin{equation}
\label{eq:banch}
[D_1, D_2]_M = [(D_1)_M, (D_2)_M].
\end{equation}
It is also clear that the derivative endomorphisms form a module over
$C^\infty(M)$. The $\mathbb R$--Lie algebra structure and the
$C^{\infty}(M)$--module structure on the vector space of derivative
endomorphisms of $E$ interact according to a Leibniz rule,
\begin{equation}
\label{eq:Leibniz}
[D_1, fD_2] = f[D_1, D_2] + (D_1)_M(f)D_2 \ ,
\end{equation}
which is a property characteristic of the sections of a
Lie algebroid. We recall the concept of Lie algebroid, introduced by
Pradines in \cite{Pradines:2}. 

\begin{df}                                  \label{df:la}
Let $M$ be a manifold. A \emph{Lie algebroid} on $M$ is a vector bundle
$(A,q,M)$ together with a vector bundle map $a\co A \to TM$ over $M$,
called the \emph{anchor} of $A$, and a bracket $[\ ,\ ]\co
\Gamma A \times \Gamma A \to \Gamma A$ which is $\mathbb{R}$--bilinear and
alternating, satisfies the Jacobi identity, and is such that
\begin{gather}
a([X,Y]) =  [a(X),a(Y)], \label{eq:la1} \\
{[X,fY]} =  f[X,Y] + a(X)(f)Y, \label{eq:la2}
\end{gather}
for all $X,Y \in \Gamma A,\ f \in C^\infty(M)$.

Given Lie algebroids $A$ and $A'$ on the same base $M$,  a \emph{Lie algebroid
morphism from $A$ to $A'$} is a vector bundle morphism $\phi\co A\to A'$
over $M$ such that 
\begin{gather}
a'\circ\phi = a, \label{anchors} \\
\phi([X, Y]) = [\phi(X), \phi(Y)],
\end{gather} 
for all $X, Y\in\Ga A$. 
\end{df}

For Lie algebroids $A$ and $A'$ on different base manifolds, the notion of
morphism is considerably more complicated; see \cite{HigginsM:1990a}. 

Any Lie algebra is a Lie algebroid over $M = \{\cdot\}$, and for any
manifold $M$, the tangent bundle $TM$ is a Lie algebroid with anchor
$a = \mbox{Id}_{TM}$. Further examples will be introduced below.

\begin{thm}\label{liealgd}
Given a vector bundle $E$ on base $M$, there is a Lie algebroid
$\D(E)$ on $M$, the smooth sections of which are the derivative 
endomorphisms of $\Ga E$, for which the bracket is the commutator 
bracket {\rm (\ref{eq:bracket})} and 
the anchor is the map $D\mapsto D_M$. 
\end{thm}

\pf
The first-order differential operators are the sections of a vector bundle
${\rm Diff}^1(E)$ on $M$. 
The symbol map is a vector bundle 
morphism, $\sigma\co{\rm Diff}^1(E)\to \Hom(T^*M, \End(E))$, 
where ${\rm End}(E)$ is the vector bundle over $M$ the fibers of which
are the endomorphisms of the fibers of $E$. 
It is a surjective submersion and has as its kernel the differential 
operators of order $0$. There is thus a short exact sequence of vector 
bundles over $M$,
$$
{\rm End}(E) \stackrel{\subseteq}{\inj} {\rm Diff}^{1}(E)
\stackrel{\sigma}{\surj}  {\rm Hom}(T^*(M), {\rm End}(E)) \ .
$$
Identifying ${\rm Hom}(T^*(M), {\rm End}(E))$ 
with $TM\otimes {\rm End}(E)$ 
in the canonical way, we define $\D(E)$ to be the 
pullback vector bundle defined by the symbol map and the injection
$TM \to TM\otimes {\rm End}(E),\ X\mapsto X\otimes{\mbox Id}_E$,
according to the diagram: 
$$
\begin{matrix}
          &\D(E)&\lrah&TM&\cr
          &&&&\cr
          &\Big\downarrow&&\Big\downarrow&\cr
          &&&&\cr
          &{\rm Diff}^1(E)&\lrah&TM\otimes{\rm End}(E).&\cr
          &&\sigma&&\cr
          \end{matrix}
$$
The pullback exists because $\sigma$ is a surjective submersion. 
Since the
right-hand vertical arrow is
an injective immersion, it follows that so is the left-hand
arrow, and we can therefore regard $\D(E)$ as a vector
subbundle of ${\rm Diff}^{1}(E)$. Similarly, because $\sigma$ is a
surjective submersion, it follows that so is the top arrow, which we
denote by $a$. Clearly the kernel of $a$ is still {\rm End}(E), and
there is an exact sequence of vector bundles over $M$,
$$
{\rm End}(E) \inj \D(E) \buildrel a \over \surj TM ,
$$
where the sections of $\D(E)$ are those first-order
differential operators $D\co \Gamma E \to \Gamma E$ for which
there exists a vector field $D_M = a(D)$ on $M$ such that
\eqref{eq:D} is satisfied.

We have already observed that for any pair of 
derivative endomorphisms $D_1,D_2$ in 
$\Ga E$, the 
bracket (\ref{eq:bracket}) is also a derivative endomorphism of $\Ga E$, with
$a([D_1,D_2]) = [a(D_1),a(D_2)]$, and that (\ref{eq:Leibniz}) is satisfied. 
Thus the proof that $\D(E)$ is a Lie algebroid on $M$ with anchor $a$ 
is complete. 
\boom

In general, a Lie algebroid is called \emph{transitive} if the anchor is 
surjective. It was shown in the course of the proof of Theorem \ref{liealgd} 
that $\D(E)$ is transitive. A right-inverse $\nabla\co TM\to\D(E)$ to the 
anchor is a linear connection in $E$ in the standard 
sense. The curvature $R_\nabla\co TM\oplus TM\to {\rm End}(E)$ of
$\nabla$ is given by 
\begin{equation}
\label{eq:curv}
R_\nabla(X_1, X_2) = \nabla_{[X_1, X_2]} - [\nabla_{X_1}, \nabla_{X_2}].
\end{equation}
In particular, the
connection is flat if and only if $\nabla$ is a morphism of Lie algebroids. 

Now consider the case of a trivial vector bundle, $E = M \times V$, and 
define a morphism from the vector bundle $TM \oplus (M \times \mathfrak{gl}(V))$
into $\D(E)$ by
$$
(X + u )(\psi ) = X(\psi ) + u(\psi )\co M \to  V \ ,
$$
where $X(\psi )$ is the Lie derivative. For an arbitrary vector
bundle, such a direct 
sum decomposition of $\D(E)$ is given locally by a trivialization of
$E$, while a global decomposition is defined by the choice of a linear
connection, $\nabla$, on $E$. Let $\nabla$ also denote the extension of the 
chosen connection to the sections of ${\rm End}(E)$. Then for $D_1 = X_1 + u_1$
and $D_2 = X_2 + u_2$, 
$$
[D_1,D_2] = [X_1,X_2] + \left( \nabla_{X_1}u_2 - \nabla_{X_2}u_1 + [u_1,u_2]
        -R_{\nabla}(X_1, X_2) \right) \ . 
$$
By definition, when $R_\nabla = 0$, this bracket gives $TM\oplus{\rm
  End}(E)$
the {\em trivial Lie algebroid structure}.
Thus the Lie algebroid $\D(E)$ is trivial if and only if $E$ admits a
flat connection. 

\medskip

The following result follows from a theorem of Bourbaki 
\cite[III.10.9, Prop. 14]{Bourbaki:algebra}.

\begin{prop}
The derivative endomorphisms 
of $\Ga E$ are the $\R$--linear
endomorphisms 
of $\Ga E$ which can be extended to derivations of the algebra of
sections of the
tensor algebra of~$E$. 
\end{prop}

The notion of derivative endomorphism is in essence purely algebraic.
Given a left module $\mathcal{E}$ over a ring 
$\mathcal{C}$, a \emph{derivation of} $\mathcal{E}$ \emph{over a
derivation $D_\mathcal{C}$ of} $\mathcal{C}$ is an additive map
$D\co \mathcal{E}\to\mathcal{E}$ such that 
\begin{equation}                                    \label{eq:Dalg}
D(f\psi) = fD(\psi) + D_\mathcal{C}(f)\psi
\end{equation}
for all $f\in\mathcal{C}$ and $\psi\in\mathcal{E}$. 
We can also refer to a derivation of $\mathcal{E}$ over a
derivation $D_\mathcal{C}$ of $\mathcal{C}$ as a \emph{derivative
endomorphism} of $\mathcal{E}$. Note, however, that in the general case, 
$D$ does not determine $D_\mathcal{C}$. 

Bourbaki \cite{Bourbaki:algebra} defines a very general notion of derivation
(involving six possibly different modules over a ring, 
three bilinear maps and three linear maps), and shows that 
the set of those derivations that belong to the ``cas I'' is closed
under the commutator. The definition of the derivations in ``cas I''
includes the derivative endomorphisms as a particular case.
The derivative endomorphisms of $\mathcal{E}$ are also a special case of the  
\emph{pseudo-linear endomorphisms}, introduced and studied in the case of a
finitely generated module over a ring by Jacobson in 
\cite{Jacobson:1935}. In fact, Jacobson uses the term `vector space 
over a field' for module over a ring and the term
\emph{transformation} for endomorphism. His definition of 
the \emph{pseudo-linear transformations} involves an 
automorphism $f \mapsto \overline f$ of $\mathcal C$ 
such that the following generalization of \eqref{eq:Dalg}
is satisfied,
$$
D(f\psi) = {\overline f} D(\psi) + D_{\mathcal{C}} (f)\psi.
$$
This generalization is applicable, for instance, when one considers 
modules over rings of complex-valued functions, in which case 
$f  \mapsto \overline f$ is the conjugation of complex numbers. 
Jacobson gave what amounts to the covariant derivation of vector fields 
in $\mathbb R^n$ as an example of such transformations, and derived the 
transformation law for connections under a change of basis.
Then he studied the generalization of the theory of elementary 
divisors to the case of these pseudo-linear transformations.
In \cite{Jacobson:1937}, he called \emph{differential transformations}
those pseudo-linear transformations for which $f = \overline f$.

Working in this algebraic setting, Herz \cite{Herz:1953}, \cite{Herz:1953II}, 
referring to Jacobson, showed that the differential transformations of a 
vector space over a commutative field form a \emph{Lie pseudo-algebra}, 
a notion which he introduced. This may now be seen as an algebraic form of 
the notion of Lie algebroid in which vector bundles over manifolds are 
replaced  by modules over rings, vector fields by derivations of
rings, and so on. 
Herz studied the sub-Lie-pseudo-algebras of a Lie pseudo-algebra over a
skew field, and then showed that what he called the kernel of the Lie 
pseudo-algebra is a Lie algebra. (There is a corresponding result 
in the theory of Lie algebroids, that the kernel of the anchor in a 
transitive Lie algebroid is a Lie algebra bundle. 
However, the algebraic results proved for Lie pseudo-algebras do not 
automatically imply the corresponding theorems in the category of Lie
algebroids, since submodules of finitely generated projective modules
do not necessarily correspond to vector subbundles.) 

A Lie pseudo-algebra that generalizes the Lie algebroid $\D(E)$ to the more 
general setting of modules over Lie pseudo-algebras was constructed in 
\cite{Mackenzie:thesis}; the construction is also given in 
\cite{Huebschmann:1990}. For references dealing with Lie pseudo-algebras, 
also now known as \emph{Lie-Rinehart algebras}, see \cite{Mackenzie:1995}. 
Rubtsov \cite{Rubtsov:1980} also introduced 
the notion of a  {\it derivation of an $A$-module over a derivation of
  an algebra $A$} in the study of the algebraic
counterpart of Lie algebroid cohomology.

While the algebraic properties of covariant derivatives were already stated 
explicitly by Schouten in 1924, in his book, {\em  Ricci Kalk{\"u}l},
the algebraization of the theory of connections is due to 
Koszul \cite{Koszul:1960}, and a very general formalization was later
accomplished by a student of Ehresmann, C.~M. de~Barros, in his thesis 
\cite{Barros:1964}. 
There, he introduced the \emph{translations
infinit{\'e}simales gradu{\'e}es}, of which the derivative
endomorphisms are a particular case, and he studied 
the \emph{repr{\'e}sentations infinit{\'e}simales gradu{\'e}es}, which 
generalize the \emph{lois de d{\'e}rivation} of Koszul.
(We shall deal with such Lie algebra representations in the
next section.) When the module being considered is that of the sections 
of a vector bundle, the derivative endomorphisms and their higher-order 
analogues were introduced by 
Palais \cite{Palais:1968} under the name \emph{quasi-scalar differential 
operators} in his analysis of elliptic self-adjoint operators on Riemannian 
vector bundles. At about the same time, Ng{\^o} Van Que \cite{NVQ:1967}, 
\cite{NVQ:1968}, \cite{NVQ:1969} used operators satisfying (\ref{eq:D}) 
in his work on Lie groupoids. 
In \cite{Kosmann-Schwarzbach:1976}, 
the name \emph{derivative endomorphisms} of $\Ga E$
for operators satisfying (\ref{eq:D}) was introduced,
and it was remarked that such operators are the
first-order differential operators on $E$ 
with scalar symbol. Such operators were independently introduced in
\cite{Mackenzie:LGLADG}, \cite{Mackenzie:1995} where, because of the 
important example \ref{example}, they are called \emph{covariant differential 
operators}. A number of papers \cite{BeilinsonS:1988}, 
\cite{ReshetikhinVW:1996} introduce an \emph{Atiyah algebra} 
of a vector bundle $E$, which is again the module of sections of $\D(E)$.
 
\medskip
 
Derivative endomorphisms arise whenever 
the infinitesimal of a representation
of a group defined by an action on a
vector bundle is considered. 
In \cite{Hermann}, Hermann introduced
the first-order differential operators 
satisfying relation \eqref{eq:D} and he briefly indicated
``the geometric genesis of these operators''.
In \cite{Kosmann-Schwarzbach:1976}, the \emph{semi-linear endomorphisms} of 
$\Ga E$ were introduced (see Section \ref{groups})
and it was shown that an action of
a Lie group on the sections of a
vector bundle by semi-linear transformations 
differentiates to an action of its Lie algebra 
by differential operators which are derivative
endomorphisms. It was stated (p.~86) that 
the Lie algebra of derivative endomorphisms of
$\Ga E$ is to be considered ``as the 
Lie algebra of the `infinite-dimensional Lie group' of all 
automorphisms of the vector bundle $E$''. 
One may say that a derivative endomorphism of $\Gamma E$ is an
\emph{infinitesimal automorphism} of the vector bundle $E$,
in the same way that a vector
field is an infinitesimal automorphism of a manifold.
Derivative endomorphisms play a crucial role throughout the theory of Lie
algebroids; it was argued in detail in \cite{Mackenzie:LGLADG} that $\D(E)$ is 
the correct generalization to vector bundles of the general linear Lie algebra 
$\mathfrak{gl}(V)$ of a vector space $V.$ 
Theorem \ref{liealgd} was proved in \cite{Mackenzie:LGLADG}.
We will discuss the relationships between these approaches further below.

\section{Derivative representations of Lie algebras 
and action Lie algebroids}\label{section2}
We again consider a vector bundle $(E,q,M)$.
Since the derivative endomorphisms 
of $\Ga E$ play the role of the  
infinitesimal automorphisms of the vector bundle $E$,
it is natural to define:

\begin{df}
\label{df:reps}
A \emph{derivative representation} of a Lie algebra $\gog$ on $E$ 
is a morphism of $\gog$ into the Lie algebra $\Gamma\D(E)$.
\end{df}

Associated to a derivative representation, $\rho$,
of $\gog$ on $E$, there is an infinitesimal action of $\gog$ on $E$,
$X \mapsto (\rho (X))_M$. 
Of course, if the base $M$ of $E$ is a point, we recover the usual notion
of linear representation. 

\begin{ex}\label{quantize}
\rm
Let us consider invariant prequantization as defined in
\cite{GuilleminS:1982}. Let $(M,\omega)$ be a
symplectic manifold, and suppose that $\omega$ is integral. Then
there is a complex line bundle $E$ on $M$, together with a
Hermitian fiber-metric, and a connection $\nabla$ in
$E$ with respect to which the metric is parallel, and such that
the curvature $R_{\nabla} $ of $\nabla$ satisfies
$$
R_\nabla = - i\omega.
$$
The prequantization action of $C^\infty(M)$ on sections of $E$ is the map,
$\delta\colon C^\infty(M)\to \Ga\D(E)$, defined by
$$
\delta(f)(\psi) = \nabla_{X_f}(\psi) + if\psi,
$$
where $f\in C^\infty(M),\ \psi\in\Ga E$, and where
$X_f$ is the Hamiltonian vector field
of $f$. It maps the Poisson bracket to the
Lie algebroid 
bracket of $\D(E)$,
respects the anchors in
the sense that $a(\delta(f)) = X_f$ for $f\in C^\infty(M)$, 
and takes values in the Lie subalgebroid $\D_0(E)$ of $\D(E)$,
whose sections respect the metric.

Now suppose that $G\times M\to M$ is a Hamiltonian action of a Lie group $G$ on
$M$, with a bracket-preserving moment map, 
$\widehat{J}\colon\gog\to C^\infty(M)$.
Then $X_{\widehat{J}(X)} = X_M$, for all $X\in\gog$, where
$X\mapsto X_M,\ \gog\to\Gamma TM$, is the infinitesimal action. Now
$\rho = \delta\circ\widehat{J}\colon\gog\to\Ga\D(E)$ is a derivative
representation of $\mathfrak g$ on $E$.
\end{ex}

There is also a well-known 
notion of a Lie algebroid representation, which we recall.

\begin{df}
\label{df:reps2}
A \emph{representation} of a Lie algebroid $A$, also on base
$M$,  on the vector bundle $E$ is a Lie algebroid morphism from $A$ 
to $\D(E)$.
\end{df}

If the base $M$ is a point,
$A$ is merely a Lie algebra and $E$ is a vector space, and this definition
reduces to that of a linear representation.
 
To express the relationship between these notions, 
we recall the definition of an \emph{action Lie algebroid}.
Consider an infinitesimal action of a Lie algebra $\mathfrak{g}$ on a 
manifold $M$, {\it i. e.}, a map  
$X\mapsto X_M$, $\mathfrak g \to \Gamma TM$, which is  
$\mathbb{R}$--linear, and
preserves brackets, $[X,Y]_M = [X_M, Y_M]$ for all
$X,Y\in\mathfrak{g}$. Extend this notation to maps $V\co M\to\mathfrak{g}$
so that $V_M$ 
is the vector field on $M$ defined by $V_M(m) = (V(m)_M)(m)$, for 
$m \in M$. Then the trivial vector bundle $M\times\mathfrak{g}$ on $M$ 
acquires a Lie algebroid structure with anchor
$a\co M\times\mathfrak{g}\to TM$ defined by $a(m,X) = X_M(m)$, and bracket
\begin{equation}                   
\label{eq:actbracket}
[V,W] = V_M(W) - W_M(V) + [V,W]^{\bullet},
\end{equation}
where $V_M(W)$ denotes the Lie derivative with respect to $V_M$
of the vector-valued function $W$, and $[V,W]^{\bullet}$ is the pointwise
bracket of maps into $\mathfrak{g}$. With this Lie algebroid structure, 
the trivial vector bundle 
$M\times\mathfrak{g}$  
is called the 
\emph{action Lie algebroid} corresponding to the infinitesimal action 
$X\mapsto X_M$, and we denote it by
$\mathfrak{g}\act M$.

Observe that if $V,W\co M\to\mathfrak{g}$ are constant maps, then the Lie
algebroid bracket $[V,W]$ is constant also, and is constant at the bracket
in $\mathfrak{g}$ of the values of $V$ and $W$. In fact this property,
together with the anchor $a\co \mathfrak{g}\act M \to TM$ and
(\ref{eq:la2}), determines (\ref{eq:actbracket}).

We can now relate 
the two concepts in Definitions \ref{df:reps} and \ref{df:reps2}. 

\begin{thm}                                
\label{thm:drep}
Let $X \mapsto X_M$ be an infinitesimal action of a Lie algebra $\mathfrak g$ 
on $M$. There is a one-to-one correspondence between derivative 
representations of $\mathfrak g$ on $E$, associated to the given 
infinitesimal action, and representations of the action Lie algebroid 
$\mathfrak{g}\act M$ on~$E$. 
\end{thm}

\pf
Let $\rho \co \mathfrak g \to \Ga \D (E)$ be a 
derivative representation of $\gog$ associated to the given infinitesimal
action $X \mapsto X_M$.
By definition, the anchor of $\rho(X)$ is $X_M$. We define a vector
bundle morphism $\sigma$ from $M \times \gog$ to $\D (E)$ by 
\begin{equation} \label{morphism}
\sigma(m,X)=(\rho(X))(m).
\end{equation}
This vector bundle morphism is a morphism of Lie algebroids from
$\gog \act M$ to $\D (E)$ because (i) it clearly satisfies the condition
on the anchors, and (ii) the bracket condition is satisfied for
constant maps from $M$ to $\gog$. 
Conversely, if $\sigma$ is a representation of $\gog \act M$ on $E$,
formula \eqref{morphism} defines a morphism $\rho$ from $\gog$
to $\Ga \D(E)$.~\boom

The significance of Theorem \ref{thm:drep} is that, on the one hand, it
can be immediately generalized once we extend the definition of 
action Lie algebroid to the case of Lie algebroids that are not Lie
algebras, and on the other hand, when the given actions
globalize, it admits a global formulation. We deal with these two aspects in
the following sections.

\begin{ex}\label{quantize2}
\rm 
If $\rho$ is the map defined in Example \ref{quantize},
the map $\sigma$ defined by \eqref{morphism}
is a representation of $\gog \act M$ on $E$, more precisely, a
morphism of Lie algebroids from $\gog \act M $ to the Lie subalgebroid
$\D_0(E)$ of $\D(E)$.
\end{ex}

\begin{ex}\rm
Given any principal bundle $P(M, H)$ and a linear representation of $H$ on 
a vector space $V$, let $E = P \times_H V$ be the associated vector bundle.
Let $A = TP/H$ be the \emph{Atiyah Lie algebroid} 
of the principal bundle. 
There is a representation $\sigma$ of $A$ on $E$ defined as 
follows (see, {\it e.g.}, \cite[App.A]{Mackenzie:LGLADG}). A section $\psi$
of $E$ can be identified with an $H$--equivariant map,
$\Tilde{\psi}\co P\to V$, and a section $X$ of $TP/H$ can be identified with 
an $H$--invariant vector field $\Tilde{X}$ on $P$. Since the Lie derivative 
of $\Tilde{\psi}$ with respect to $\Tilde{X}$ is also an $H$--equivariant 
map from $P$ to $V$, it corresponds to a section $\sigma(X)(\psi)$. 

If, in particular, 
$P$ is a Lie group $G$ with $H$ a closed subgroup, this construction
associates to any linear representation of $H$ on $V$ the
representation $\sigma$ of $TG/H$ on $E =G\times_H V$. On the
other hand, the Lie algebroid $TG/H$ is canonically isomorphic to
the action Lie algebroid $\mathfrak{g}\act (G/H)$ arising from the standard action
of $\mathfrak{g}$ on $M = G/H$, under the map 
$TG/H \to \mathfrak{g} \act (G/H)$ defined by the canonical right
identificatiion of $TG$ with $G \times \mathfrak{g}$. 
It follows by Theorem \ref{thm:drep} that there is a one-to-one
correspondence between the linear representations of $\mathfrak{h}$ 
on $V$ and the derivative representations of $\gog$ on $E = G \times_H  V$. 

In this case, the derivative representation of $\mathfrak g$ on 
$ E= G \times_H V$
associated to the representation $\sigma$ of $A=TG/H$ on $E$
is the infinitesimal of the linear 
representation of $G$ induced from that of $H$.
So this result is the infinitesimal 
counterpart of an observation 
that was made in
\cite{Kosmann-Schwarzbach:1976}
relating semi-linear representations of Lie groups to induced representations.
\end{ex}

\begin{ex}\rm
A representation of a Lie algebroid $A$ on a vector bundle 
$E$ extends to a representation by derivations of its exterior algebra, that
is, a morphism from the $\R$--Lie algebra $\Ga A$ to the $\R$--Lie algebra
of derivations of $\Ga(\bigwedge E)$, and conversely the restriction
to $\Ga E$ of such a representation by derivations of $\Ga(\bigwedge E)$
is a representaion of $A$ on $E$ in the sense of Definition \ref{df:reps2}.
\end{ex}

\begin{rmk} \rm
Consider a Lie algebroid $A$ and a vector bundle $E$ on the same base $M$. 
Several authors \cite{Xu:1999} \cite{EvensLW} consider the
$A$--\emph{connections on} $E$, defined as the
vector bundle morphisms $A \to \D(E)$ satisfying  
condition \eqref{anchors} on the anchors, without a bracket condition. 
If $A = TM$, this is precisely the standard concept of linear connection
in $E$. 

The \emph{curvature} of such an $A$--connection $\nabla$ is defined as in 
(\ref{eq:curv}). Now a representation  of $A$ on $E$ is precisely a 
flat $A$--connection on $E$.  Extensive use of the notion of $A$--connection 
and of representation of Lie algebroids was made in \cite{Xu:1999} 
\cite{EvensLW}, especially the representations of $A$ on $\bigwedge^n A$ 
where $n$ is the rank of $A$. The one-to-one correspondence described in 
the proof of Theorem \ref{thm:drep} is actually the restriction to the 
derivative endomorphisms of a one-to-one correspondence between
derivative pre-representations of $\gog $ on $E$, where we only
require linearity of the map $\gog \to \Ga \D(E)$, and $A$--connections
on $E$. 
\end{rmk}

\medskip

The vector bundle that 
we have called the Atiyah Lie algebroid of a principal bundle 
was introduced by Atiyah in 1957 \cite{Atiyah:1957}, ten years before 
Pradines introduced the abstract concept of Lie algebroid. It is also 
called the \emph{gauge algebroid} of the principal bundle, and its sections 
are then called \emph{infinitesimal gauge transformations} of the bundle. 
The definition of derivative representation comes from 
\cite{Kosmann-Schwarzbach:1976}. 
The derivative representations are a particular case of the graded
infinitesimal representations of de Barros \cite{Barros:1964}.
When $\mathfrak g$ is the Lie algebra of derivations of
$C^{\infty}(M)$
and when the map $\mathfrak g \to \Gamma {\mathfrak D} (E)$ is 
$C^{\infty}(M)$--linear, a derivative representation is a flat
derivation law in the sense of Koszul \cite{Koszul:1960}. 
The definition of Lie algebroid representation 
may be found in \cite[III.2.9]{Mackenzie:LGLADG}; it was probably
first formulated in \cite{NVQ:1969}. The corresponding algebraic
notion is that of a left module over a Lie pseudo-algebra.

\section{Lie algebroid actions and representations} \label{sect:drLa}

In this section we first consider the infinitesimal actions of a Lie 
algebroid $A$ with base $M$ on a fibered manifold with base $M$, that is, 
on a surjective submersion onto $M$. 
Such an infinitesimal action of $A\to M$ on $\varphi\co F \to M$ defines a Lie
algebroid structure on the vector bundle $\varphi^* A$, the pullback
of $A$, which is called the action Lie algebroid
defined by the given infinitesimal action of $A$ on $F$.

We then extend the definition of the derivative representations
to the case of a Lie algebroid $A$ and a vector bundle with base $F$.
As a particular case of this notion, we recover Definition \ref{df:reps}
in the case where $M$ is a point. Finally in this section we show that 
the one-to-one correspondence between derivative representations and 
representations of action Lie algebroids can be extended to this more 
general situation.

\begin{df}                                       
\label{df:algdaction}
Let $A$ be a Lie algebroid on $M$, and let $\varphi\co F\to M$ be a
fibered manifold. An \emph{infinitesimal action of $A$ on $F$} is an
$\R$--linear map $X\mapsto X_F,\ \Ga A\to \Ga TF$, such that

\noindent
{\rm (i)} for each $X \in \Ga A$, $X_F$ is projectable to $a(X)$,

\noindent
{\rm (ii)} the map $X \mapsto X_F$ preserves brackets,

\noindent
{\rm (iii)} the map $X \mapsto X_F$ is $C^{\infty}(M)$--linear in the following
sense: for each $f \in C^{\infty}(M)$ and each $X \in \Ga A$,
\begin{equation}\label{module}
 (fX)_F = (f\circ \varphi) X_F.
\end{equation}
\end{df}

If $\Gamma TF$ is given the $C^{\infty}(M)$--module structure defined by 
$(f,Y) \mapsto (f \circ\varphi) Y$, $C^{\infty}(M) \times \Gamma TF
\to \Gamma TF$, then Equation \eqref{module} is the condition that 
$X \mapsto X_F$ is a morphism of left modules.

If, in particular, the base $M$ is a point, then $A$ is a Lie algebra, and 
we recover the notion of infinitesimal action of a Lie algebra on a manifold $F$,
conditions (i) and (iii) being trivially satisfied.

If $F=M$ and $\varphi$ is the identity, an infinitesimal action of 
the Lie algebroid $A$ on the manifold $M$
is a $C^{\infty}(M)$--linear map $X \mapsto X_M$ from $\Ga A$ to
$\Gamma TM$, 
{\it i.e.}, a morphism of vector bundles from $A$ to $TM$, which by
condition (i) is necessarily the anchor of $A$.

It is shown in \cite{HigginsM:1990a} that whenever an infinitesimal 
action of $A$ on $\varphi\co F\to~M$ is defined, there is an associated 
Lie algebroid structure on the pullback vector bundle, $ \varphi^* A$,
$$
\begin{matrix}
          &\varphi^*A&\lrah&A&\cr
          &&&&\cr
          &\Big\downarrow&&\Big\downarrow&\cr
          &&&&\cr
          &F&\lrah&M.&\cr
          &&\varphi&&\cr
          \end{matrix}
$$
With this Lie algebroid structure, $\varphi^* A$ is called the \emph{action
Lie algebroid} associated to the infinitesimal action $X \mapsto X_F$, 
and it is denoted $A \act F$. We recall the definition of this Lie algebroid 
structure. We identify the module of sections of the pullback with
$C^\infty(F)\otimes\Ga A$, where the tensor product is over $C^\infty(M)$. 
Thus $h\otimes fX = h(f\circ \phi)\otimes X$ and $h\otimes X$ is 
identified with $ h (X\circ \varphi)$,
for $h\in C^\infty(F),\ f\in C^\infty(M), \ X\in \Ga A$. 
Then $\varphi^* A$ is a Lie algebroid with the following anchor, $a_F$, 
\begin{equation}
a_F\left ( h\otimes X\right)  =   h X_F \ ,
\end{equation}
and $\mathbb R$--bilinear bracket,
\begin{equation}                   \label{eq:actionbracket}
\left[ h\otimes X, k\otimes Y\right] = 
h k\otimes {[}X,Y{]} 
+  h X_F(k)\otimes Y - k Y_F(h)\otimes X, 
\end{equation}
where $h, k \in C^\infty(F),\ X,Y\in\Ga A$. 

If the base $M$ is a point, this Lie algebroid over $F$ is merely the
product $F \times A$ with the action Lie algebroid structure over base
$F$ defined in Section \ref{section2}.

The notion of action and the construction of an action Lie algebroid
extend to the situation in which $\phi$ is an arbitrary smooth map, 
but we shall not need this generality here.

\begin{ex}\rm
Let $A$ be any Lie algebroid on a connected base $M$, and let
$f\co\Tilde{M}\to M$ be any covering, with group $\pi$. Then
there is a canonical infinitesimal action of $A$ on $\Tilde{M}$ 
in which $X_{\Tilde{M}}$
is the $\pi$--invariant lift of $a(X)$ to $\Tilde{M}$.
\end{ex}

\begin{ex}\rm
Given any principal bundle $P(M,H)$, the Atiyah Lie algebroid
$TP/H$ with base $M$ acts on $P\to M$ by
$X\mapsto\Ri{X}$, where $\Ri{X}$ is the right-invariant vector field on
$P$ defined by $X\in \Ga (TP/H)$; 
the action Lie algebroid is canonically isomorphic to $TP$.
\end{ex}

In view of Theorem \ref{thm:drep}, we shall define the notion of 
derivative representation of a Lie algebroid $A$ with base $M$
on a vector bundle $E$ with base $F$, associated to an infinitesimal 
action of $A$ on $F$, in such a way that the one-to-one correspondence 
between derivative representations and representations of the action Lie 
algebroid is extended.

\begin{df}
\label{df:derep}
Let $A$ be a Lie algebroid on $M$, and let $\varphi \co F \to M$ be a
fibered manifold. Let $X \mapsto X_F,\ \Ga A \to \Gamma TF$,
be an infinitesimal action of $A$ on $F$. Let us consider a vector 
bundle $E$ on $F$. A \emph{derivative representation of $A$ on $E$} associated 
to the infinitesimal action $X \mapsto X_F$ is a morphism, $\rho$, of Lie
algebras from $\Ga A$ to $\Ga \D (E)$ such that 

\noindent
{\rm (i)} 
for any $X \in \Ga A$,
\begin{equation} \label{associated}
(\rho(X))_F = X_F \ ,
\end{equation}

\noindent
{\rm (ii)} $\rho$ is $C^{\infty}(M)$--linear
in the following sense: for each $f \in C^{\infty}(M)$ and each $X\in
\Ga A$, 
\begin{equation}\label{eq:third}
\rho(fX)= (f \circ \varphi) \rho(X) \ .
\end{equation}
\end{df}

Since for any $X \in \Ga A$, $\rho(X)$ is a section of $\D(E)$, 
condition \eqref{associated} expresses the fact that, 
for any $h\in C^\infty(F),\psi\in\Ga E$,
\begin{equation}\label{eq:first}
\rho(X)(h\psi) = h \rho(X)(\psi) + X_F(h)\psi.
\end{equation}

\begin{thm}
\label{thm:drepoid}
Let $A$, $F$, $X \mapsto X_F$ and $E$ be as in Definition \emph{\ref{df:derep}}.
There is a one-to-one correspondence between derivative
representations of $A$ on $E$ associated to the infinitesimal action
$X \mapsto X_F$ and representations of the action Lie algebroid 
$A \act F$ on $E$.
\end{thm}

\pf
We sketch the proof which is very similar to that of Theorem \ref{thm:drep}. 
If $\rho: \Ga A \to \Ga \D(E)$ is a derivative representation of $A$
on $E$ associated to the given infinitesimal action of $A$ on $F$,
we define $\sigma \co \varphi^*A \to E$ in terms of sections of the
form $h \otimes X$, $h \in C^{\infty}(F), X \in \Ga A$, by
$$
\sigma (h \otimes X) (p)= h(p) \rho(X)(p) \ , 
$$
for any $p \in F$. Then $\sigma$ is a representation of the Lie 
algebroid $\varphi^*A$ on $E$.

Conversely, if $\sigma$ is a representation of the Lie algebroid $\varphi^*A$ on
$E$, we set, for $X \in \Ga A$, and $p \in F$,
$$
\rho(X) (p)= \sigma(X \circ \varphi)(p) \ .
$$
Then $\rho$ is a derivative representation of $A$ on $E$.
\boom

Thus we have obtained a generalization of Theorem \ref{thm:drep} to the 
case of Lie algebroids acting on vector bundles.

The Lie algebroid structure \eqref{eq:actionbracket} of an action Lie 
algebroid corresponds in the algebraic setting treated by Fel'dman 
\cite{Feldman:1982} to a \emph{crossed product} structure on its 
space of sections.

\section{Semi-linear representations of Lie groups and action groupoids}
\label{groups}
In this section we provide a global formulation in terms of 
Lie groupoids for the construction described in Section \ref{section2}:
the relationship between derivative representations of Lie algebras 
on vector bundles and action Lie algebroids has a global analogue,
the relationship between the semi-linear representations of Lie groups 
acting on vector bundles and the action groupoids, described in Theorem
\ref{thm:slrep} below. In the following section we shall discuss the more 
general case of a global formulation for Lie algebroid representations.
These global constructions are such that applying the Lie functor to the
Lie groups and Lie groupoids under consideration yields the known 
infinitesimal constructions for Lie algebras and Lie algebroids.

For the basic definitions of groupoid theory used in what follows, see 
\cite{Mackenzie:LGLADG} or \cite{Landsman:1998} or \cite{CannasdaSilvaW} 
(note that the composition convention in 
\cite{CannasdaSilvaW} is the opposite 
to that followed here). Given 
a vector bundle $(E,q,M)$, we denote by $\Phi(E)$ the set of all 
linear isomorphisms from a fiber of $E$ to some (in general different) 
fiber of $E$. Then $\Phi(E)$ is a Lie groupoid with base $M$, 
with its groupoid structure 
arising from the ordinary composition of maps, and its smooth structure 
given by locally identifying $\Phi(M\times V)$ with $M\times GL(V)\times M$. 
(In \cite{Mackenzie:LGLADG}, this groupoid is called the \emph{linear 
frame groupoid} of $E$.) The theorem which follows is fundamental to the
results of this section. 

\begin{thm}\label{LAofPhi}
The Lie algebroid of the Lie groupoid $\Phi(E)$ is canonically
isomorphic to $\D(E)$.
\end{thm}

Theorem \ref{LAofPhi} emerged from two separate lines of work. 
On the one hand, a proof was given by Kumpera \cite{Kumpera:1971} in
unpublished notes of 1971, but was not included in the published
version (Appendix A of \cite{KumperaS}). A simplification of Kumpera's 
proof was given by Mackenzie in 1987 \cite[III.4.5]{Mackenzie:LGLADG}.
On the other hand, it was remarked by Hermann in \cite{Hermann}
and independently by Kosmann-Schwarzbach
in \cite{Kosmann-Schwarzbach:1976} that the derivative endomorphisms
of the space of sections of
a vector bundle are the infinitesimal generators of one-parameter groups 
of vector bundle automorphisms. This leads precisely to the result 
that $A\Phi(E)\isom\D(E)$. 

\medskip

We now define the semi-linear isomorphisms of $\Gamma E$ and
correspondingly the semi-linear representations of a Lie group on the
vector bundle, $E$.

\begin{df}\label{sl}
Let $(E,q,M)$ be a vector bundle. A \emph{semi-linear isomorphism} 
of $\Ga E$ is an $\R$--linear automorphism $\mu$ of 
$\Ga E$ such that there exists an $\R$--linear isomorphism, $\mu^M$,
of $C^\infty(M)$ satisfying 
\begin{equation}\label{semilinear}
\mu(f\psi)= \mu^M (f) ~ \mu(\psi) \ ,
\end{equation}
for all $f \in C^{\infty}(M)$ and all $\psi \in \Gamma E.$
\end{df}

It follows that $\mu^M$ is a ring automorphism of $C^{\infty}(M)$.
Any automorphism $\nu$ of the vector bundle $E$ that projects onto a 
transformation $\nu_M$ of the base manifold gives rise to a semi-linear
isomorphism $\psi \mapsto \nu \cdot \psi$ of $\Gamma E$ defined by
\begin{equation}\label{isom}
\nu \cdot \psi = \nu \circ \psi \circ (\nu_M)^{-1} \ ,
\end{equation}
in which case the ring endomorphism of $C^{\infty}(M)$ is $f \mapsto f
\circ (\nu_M)^{-1}$, 
and conversely, every semi-linear isomorphism 
of $\Gamma E$ arises from a vector bundle automorphism of $E$.
(See \cite{Kosmann-Schwarzbach:1972} \cite{Kosmann-Schwarzbach:1976}.)
\begin{prop}\label{diffeo}
Formula \eqref{isom} establishes a one-to-one correspondence between
vector bundle automorphisms of $E$ and semi-linear isomorphisms of
$\Gamma E$.
\end{prop}

The semi-linear isomorphisms of $\Gamma E$ are the $\mathbb
R$--isomorphisms of $\Gamma E$ which can be extended to automorphisms of
the algebra of sections of the tensor algebra of~$E$.

\begin{df}\label{slrep}
A \emph{semi-linear representation} of a Lie group $G$ on $E$ is  
a morphism, $R$, from $G$ into the group of semi-linear isomorphisms of
$\Gamma E$ such that, for any $\psi \in \Gamma E$, the map 
$(g,m) \mapsto R(g)(\psi)(m)$, $G \times M \to E$, is smooth.
\end{df}

Each semi-linear representation, $R$, of $G$ on $E$ is associated to an
action $g \mapsto g_M$ of $G$ on $M$ such that $(R(g))^M(f)=f \circ
(g_M)^{-1}$.

There is also a well-known notion of a Lie groupoid representation, 
which we recall.

\begin{df}
\label{df:reps3}
A \emph{representation} of a Lie groupoid $G$ with base
$M$ on the vector bundle $E$ is a base-preserving Lie groupoid 
morphism from $G$ to $\Phi(E)$.
\end{df}

If the base $M$ is a point, $G$ is merely a Lie group and $E$ is a 
vector space, and this definition reduces to that of a linear representation.
In general it follows from Theorem \ref{LAofPhi} that applying the Lie functor 
to a representation of $G$ on $E$ yields a representation of $AG$ on $E$.

To express the relationship between the notion of semi-linear 
representation and that of a groupoid representation,
we recall the definition of \emph{action groupoid}.
Consider an action $g\mapsto g_M$ of a Lie group $G$ on a 
manifold $M$. The product manifold $G \times M$ becomes a Lie groupoid
with base $M$ when it is equipped with the source map $(g,m) \mapsto m$,
the target map $(g,m) \mapsto g_M(m)$, for $g \in G$, $m\in M$, the 
partial multiplication 
$((g_1,m_1),(g_2,m_2)) \mapsto (g_1 g_2, m_2)$ defined if and only if 
$m_1 = (g_2)_M(m_2)$, and the inversion $(g,m) \mapsto 
(g^{-1}, g_M(m))$. With this Lie groupoid structure, the product 
manifold $G \times M$ is called the \emph{action groupoid} 
corresponding to the action 
$g \mapsto g_M$, and we denote it by $G \act M$. The following result
is proved in \cite[Th. 2.5]{HigginsM:1990a}. 

\begin{prop}\label{LAAG}
If $\mathfrak g$ is the Lie algebra of the Lie group $G$, the Lie
algebroid of the action groupoid $G \act M$ is the action Lie
algebroid $\mathfrak g \act M$ associated to the infinitesimal of the 
action of $G$.
\end{prop}

\begin{ex}\label{quantize3} 
\rm The representation $\sigma$ of the Lie algebroid $\gog \act M$ on $E$
of Example \ref{quantize2} can be integrated to a morphism 
$G\act M\to\Phi_0(E)$, where $\Phi_0(E)$ is the Lie subgroupoid of $\Phi(E)$ of
isometries between the fibers of $E$, {\it i.e.}, to an orthogonal representation
of $G \act M$ on $E$, if and only if $(E,\nabla)$ together with the metric 
constitute $G$--invariant data in the sense of \cite{GuilleminS:1982}.
\end{ex}

We can now relate the two concepts defined in \ref{slrep} and 
\ref{df:reps3}.

\begin{thm}                                
\label{thm:slrep}
Let $g \mapsto g_M$ be an action of a Lie group $G$ on $M$.
There is a one-to-one correspondence between semi-linear representations of 
$G$ on $E$, associated to the given action, and representations of the action 
groupoid $G \act M$ on~$E$. 
\end{thm}

\pf The proof follows from the observations that semi-linear
representations are in one-to-one correspondence with the 
actions of $G$ on $E$ by vector bundle automorphisms, which in turn are in
one-to-one correspondence with morphisms of Lie groupoids from the
action groupoid, $G \act M$, to $\Phi(E)$, {\it i.e.}, representations 
of the action groupoid on $E$ \cite{Mackenzie:LGLADG}.
\boom

In view of the properties of the Lie functor stated in Theorem
\ref{LAofPhi} and Proposition \ref{LAAG}, this theorem implies 
Theorem \ref{thm:drep} when the action of $\gog$ on $M$ can be 
integrated globally. 

\begin{rmk}\rm
The semi-linear transformations of $\Gamma E$ are another particular
case of the pseudo-linear transformations of Jacobson
\cite{Jacobson:1937}, in which the linear mapping $D_{\mathcal C}$ 
vanishes, and the ring automorphism $f \mapsto \bar f$ of
$\mathcal C=C^{\infty}(M)$ is, in the notation of Definition
\ref{sl}, the map $\mu^M$.
\end{rmk}

\begin{rmk}\rm
The isomorphism of Theorem \ref{LAofPhi} arises as an infinitesimal 
version of the correspondence between vector bundle automorphisms and 
semi-linear isomorphisms (see 
Proposition \ref{diffeo}), but also 
admits another formulation. 

A vector field $X$ on a vector bundle $(E, q, M)$ is \emph{linear}
if the flow of $X$ consists of local vector bundle automorphisms of $E$. 
Such a vector field clearly projects to a vector field $X_M$ on $M$. 
In \cite{MackenzieX:1998}, linear vector fields are characterized as 
pairs $X\in\Gamma TE,\ X_M\in\Gamma TM$ which constitute a vector bundle 
morphism from $(E, q, M)$ to $(TE, T(q), TM)$. 

The one-to-one correspondence between linear vector fields on $E$
and derivative endomorphisms of $\Gamma E$  
can be expressed directly as follows \cite{Kosmann-Schwarzbach:1978}
\cite{MackenzieX:1998}.
Given a vector field $X$ on $E$ which projects to a vector field $X_M$ 
on $M$, and a section $\psi$ of $E$, we define 
\begin{equation}\label{lieder}
D_X(\psi) (m) = T_m(\psi)(X_M(m)) - X(\psi(m)),  
\end{equation}
for $m \in M$. (Here the right-hand side 
is a vertical tangent vector, which we
regard as an element of $E$ in the usual way.)
When $X$ is linear, the ${\mathbb R}$--linear map $D_X$ is a
derivative endomorphism of $\Gamma E$ associated with $X_M$, and is called the 
\emph{Lie derivation} with respect to $X$. 
Conversely, given a derivative endomorphism $D$ of $\Gamma E$
associated with a derivation $X_M$ of $C^{\infty}(M)$, we set
\begin{equation}\label{inverselieder}
X_D (p)=  T_m(\psi)(X_M(m)) - (D\psi)(m),  
\end{equation}
where $p \in E$ and $\psi$ is a section of $E$ such that $\psi(m) = p$.
Then $X_D$ is a well-defined linear vector field on $E$, and the maps
$X \mapsto D_X$ and $D \mapsto X_D$ are obviously inverses of one
another. Moreover this one-to-one correspondence preserves brackets.
In fact, the linear vector fields are the section of a Lie algebroid
$T^{{\rm LIN}}E$
with base $M$, defined in \cite{MackenzieX:1998}, and this correspondence 
arises from a Lie algebroid isomorphism
from $T^{{\rm LIN}}E$ to $\D(E)$. Therefore
the isomorphism of Theorem \ref{LAofPhi} may also be formulated as 
an isomorphism between the Lie algebroid of $\Phi(E)$ and $T^{\rm LIN}E$. 

An infinitesimal action of a Lie algebroid $A$ 
on a vector bundle $E$ in the
sense of Definition \ref{df:algdaction}
is called \emph{linear} if the action is by linear vector fields. 
Since a Lie algebroid morphism from $A$ 
to $T^{\rm LIN}E$ can be identified with a Lie algebroid 
morphism from $A$ to $\D(E)$, it follows that the linear infinitesimal
actions of $A$ on $E$ can be identified with the representations of
$A$ on $E$ in the sense of Definition \ref{df:reps2}.
Since, by condition (i) of Definition \ref{df:algdaction},
for each $X \in \Gamma A$, $X_M= a(X)$, the correspondence established
by \eqref{lieder} coincides with that of \cite[p. 212]{HigginsM:1990a}.

A more general non-linear version of the bijective correspondence
between vector fields on $E$ and differential operators on the
sections of $E$ 
was established 
in \cite{Kosmann-Schwarzbach:1980}. Consider a fiber
bundle $(F,q,M)$. Then, given an infinitesimal automorphism of the bundle, 
that is, a vector field $X$ on $F$ which projects to a vector field $X_M$ 
on $M$, and a section $\psi$ of $F$, formula \eqref{lieder} defines a 
vertical tangent vector to $F$, and we thus obtain a section
of $\psi^*VF$, where $VF$ is the vertical tangent bundle of $F$. 
In this way one obtains an operator which is a 
\emph{first-order differential section operator} from $F$ to $VF$ 
in the terminology of \cite{Palais:1968}. Conversely, 
such an operator corresponds to
an infinitesimal automorphism if and only if it is quasi-scalar. 
For a fibered manifold $(F, q, M)$, the set of infinitesimal automorphisms
has both a module structure over $C^\infty(M)$ and a bracket, making it a Lie
pseudo-algebra, but in general there is no underlying Lie algebroid. 
\end{rmk}

\section{Groupoid actions and representations}
\label{sect:gar}

More generally, we now consider the action of a Lie groupoid on a
fibered manifold. (We shall not consider 
the still more general case of the action of a Lie groupoid on an 
arbitrary smooth map. For this, see both \cite{HigginsM:1990a} and 
the recent account of Moerdijk and Mr\v cun \cite{MoerdijkM}.)

Let $G$ be a Lie groupoid with base $M$ and source map $\alpha$,
and let $\phi:F \to~M$ be a fibered manifold. We denote by $\phi^*G$ 
the pullback in the following commutative diagram,
$$
\begin{matrix}
          &\varphi^*G&\lrah&G&&\cr
          &&&&\cr
          &\Big\downarrow&&\Big\downarrow&\alpha&\cr
          &&&&\cr
          &F&\lrah&M.&&\cr
          &&\varphi&&&\cr
          \end{matrix}
$$
For $m \in M$, we denote $\phi^{-1}(m) \subseteq F$ by $F_m$.

\begin{df}\label{gpdaction}
Let $G$ be a Lie groupoid 
with base $M$, source map $\alpha$ and target map $\beta$,
and let $\phi:F \to M$ be a fibered manifold.
An \emph{action} of $G$ on $F$ is a map, $S$, from $G$ to the set of smooth
bijective maps from a fiber of $F$ to some (in general different)
fiber of $F$, such that 

\noindent
{\rm (i)} for any $g \in G$, $S(g)$ is a map from 
$F_{\alpha(g)}$ to $F_{\beta(g)}$, 

\noindent
{\rm (ii)} for any $g_1 \in G$ and $g_2 \in G$ such that $\beta(g_1)=\alpha(g_2)$, 
$S(g_2 g_1)=S(g_2) \circ S(g_1)$,

\noindent
{\rm (iii)}
the map $(g,p)\mapsto S(g)(p)$ from $\phi^*G$ to $F$ is smooth.
\end{df}

If, in particular, the base $M$ is a point, then $G$ is a Lie group,
and we recover the notion of action of a Lie group on a manifold $F$.

When $F=M$ and $\phi$ is the identity, the only action of the Lie groupoid 
$G$ on $M$ is that given by $S(g)(\alpha g) = \beta g$. 

It was originally shown by Ehresmann \cite{Ehresmann:1957} that whenever 
an action, $S$, of $G$ on $\phi:F \to M$ is defined, there is an associated 
Lie groupoid structure on the manifold $\phi^*G$, with base $F$. 
The source map is $(g,p) \mapsto p$, while the target map is 
$(g,p) \mapsto S(g)(p)$. The partial multiplication is 
$((g_1,p_1),(g_2,p_2)) \mapsto (g_1g_2, p_2)$ defined if and only if 
$g_1=S(g_2)p_2$, and the inversion is $(g,p) \mapsto (g^{-1}, S(g)(p))$.
With this Lie groupoid structure, $\phi^* G$ is called the 
\emph{action groupoid} associated to the given action, and we denote it by 
$G \act F$. The top arrow in the preceding diagram is a morphism of Lie 
groupoids from $G \act F$ to $G$, over $\phi : F \to M$.

If $M$ is a point, then the Lie groupoid $G$ is a Lie group, and the Lie 
groupoid structure defined on $G \times F$ reduces to the action groupoid 
structure described in Section \ref{groups} associated to a Lie group action 
on a manifold. 

The Lie functor maps action groupoids to action Lie algebroids. More
precisely, the differential, $\sigma$, of $S$ is defined by
$$
\sigma(X)(p) = TS(p)(X_{\phi(p)}) \ ,
$$
where $p \in F$ and 
$X \in \Gamma AG $, while $TS(p)$ denotes the tangent map at
the unity, $\phi(p) \in M$, of the map $g \mapsto S(g)(p)$ from $G$ to $F$.
The result which follows was proved in \cite{HigginsM:1990a}. 

\begin{prop}\label{diff}
Let $S$ be an action of the Lie groupoid $G$ on $\phi:F \to M$. Then
$S$ differentiates to an infinitesimal action of the Lie algebroid
$AG$ on $F$, and the Lie algebroid of the action groupoid $G \act F$ is
the action Lie algebroid $AG \act F$.
\end{prop}
 
Before considering a notion of semi-linear representation for Lie groupoids,  
we recall the notion of bisection. 

\begin{df}
Let $G$ be a Lie groupoid, with base $M$, source map $\alpha$ and target 
map $\beta$. A map $b\co M \to G$ is a \emph{bisection} of $G$ if 
$\alpha \circ b = \rm{Id}_M$ and $\beta \circ b$ is a diffeomorphism of $M$.
\end{df}

The set of bisections of $G$ is a group under the composition law
defined by
$$
(b_1 \cdot b_2)(m) = b_1((\beta \circ b_2)(m))b_2(m) \ , 
$$ 
for bisections $b_1$ and $b_2$, and $m \in M$.
The \emph{group of bisections} of $G$ is denoted by ${\mathcal B}(G)$.

\begin{prop}\label{canon}
Let $(E,q,N)$ be a vector bundle. There is a canonical group 
isomorphism from ${\mathcal B}(\Phi(E))$ to the group of vector 
bundle automorphisms of $E$. 
\end{prop}

\pf
Given a bisection, $b$, of $\Phi(E)$, we define a bundle map, 
$\tilde b$, from $E$ to itself by
$$
\tilde b(v)=b(q(v))(v) \ ,
$$
for $v \in E$. Then $\tilde b$ is a vector bundle automorphism of $E$, 
projecting onto the diffeomorphism $\beta \circ b$ of $N$. Conversely, 
given a vector bundle automorphism $u$ of $E$, the map $x \mapsto u|_{E_x}$, 
$x \in N$, defines a bisection of $\Phi(E)$. The smoothness of the map is
proved by using the local triviality of $E$. These maps are mutually inverse 
and they are group morphisms.
\boom

It follows from Proposition \ref{diffeo}
that ${\mathcal B}(\Phi(E))$ is isomorphic to the group of
semi-linear isomorphisms of $\Ga E$. 

Now consider Lie groupoids $G_1$ and $G_2$ on the same base $M$. It is clear 
that any Lie groupoid morphism $u\co G_1\to G_2$ over $M$ defines a group 
morphism $U\co{\mathcal B}(G_1) \to {\mathcal B}(G_2)$ of the corresponding 
groups of bisections by 
\begin{equation}\label{eq:uU}
U(b)(m) = u(b(m)) ,
\end{equation}
for $b \in {\mathcal B}(G_1)$ and $m \in M$. 
Thus a semi-linear representation
of a Lie group $G$ on a vector bundle $(E, q, M)$ may be regarded as a group
morphism from $\mathcal{B}(G)$ to ${\mathcal B}(\Phi(E))$ which satisfies the
smoothness condition in \ref{slrep}. 

\medskip

Let $G$ be a Lie groupoid acting on a fibered manifold,
$\phi : F \to M$, with action groupoid $G\act F$. 
We consider a vector bundle $E$ on 
base $F$ and a representation of $G\act F$ on $E$,
{\it i.e.}, a morphism of Lie groupoids from $G\act F$ to $\Phi(E)$. 
For groupoids in the category of sets, this notion was
studied by Brown \cite{Brown:1970}. In the smooth case,
such a representation is 
a global form of a derivative representation. Indeed, by 
Proposition \ref{diff} it 
differentiates to a representation of $AG \act F$ on
$E$, which by Theorem \ref{thm:drepoid} is a derivative representation of
$AG$ on $E$. 
 
Thus for each $g\in G$ and $p\in F_{\alpha g}$, we are given an 
isomorphism from $E_p$ to $E_{S(g)p}$, where $S$ is the given action
of $G$ on $F$. 
According to \eqref{eq:uU} and to Proposition \ref{canon}, there is a group 
morphism from $\mathcal{B}(G\act F)$ to the group of semi-linear isomorphisms 
of $\Gamma E$. In addition, given a bisection $b$ of $G$, there is a 
bisection $\hat{b}$ of $G\act F$ defined by 
$$
\hat{b}(p) = (b(\phi(p)), p), 
$$
for $p\in F$, and the map $b\mapsto\hat{b}$ is a group morphism from 
$\mathcal{B}(G)$ to $\mathcal{B}(G\act F)$. Composing these maps we obtain 
a group morphism from $\mathcal{B}(G)$ to the group of semi-linear isomorphisms 
of $\Gamma E$. This group morphism, $R$, is associated with the given action 
of $G$ on $F$, in the sense that if $v \in E_p$ for $p\in F$, 
then $R(b)v \in E_{p'}$, where $p'=S(b(\phi p))p$, for all 
$b \in \mathcal{B}(G)$. Thus we have:

\begin{prop}\label{prop:gpdslr}
Given an action of a Lie groupoid $G$ on a fibered manifold 
$\phi:F \to M$, to any representation of $G \act F$ on a vector bundle $E$
with base $F$ there corresponds a representation of the group of
bisections of $G$ by semi-linear isomorphisms of $\Gamma E$, 
associated to the given action.
\end{prop}

Since the Lie algebra of sections of $AG$ can be formally regarded as the 
Lie algebra of the infinite-dimensional group ${\mathcal B}(G)$, this 
proposition can be viewed as a global form of part of Proposition 
\ref{thm:drepoid}. When the Lie groupoid $G$ is a Lie group, its group 
of bisections is the group itself, so this proposition is a partial
generalization of Theorem \ref{thm:slrep} to the case of Lie groupoids. We 
may therefore regard the group morphisms from $\mathcal{B}(G)$ to the group of
semi-linear isomorphisms of $\Gamma E$ satisfying a suitable smoothness 
condition as the \emph{semi-linear representations} of $G$ on $E$. However, 
the characterization of smoothness for such semi-linear representations of 
groupoids in terms of the bisections of $G$ alone appears to be a difficult 
matter. The concept of bisection may readily be localized, leading to a sheaf of 
germs of local bisections, but the relationship between bisections of $G$ and
bisections of $G\act F$ evades a simple description. 

\begin{rmk}\rm
The representation of the group of bisections in Proposition \ref{prop:gpdslr}
may be approached in an alternative way. We consider a
representation of $G \act F$ on $E$ as above and, for $g \in G$ and
$p \in F_{\alpha g}$, we denote 
the isomorphism from $E_p$ to $E_{S(g)p}$ by $\rho(g,p)$. 
Given $g\in G$, to any section $\psi$ of $E$ 
over $F_{\alpha g}$ we can associate a section $g \cdot \psi$ of $E$ 
over $F_{\beta g}$, by
\begin{equation}\label{Rrhog}
(g \cdot \psi) (p) = \rho(g, S(g^{-1})p)(\psi(S(g^{-1})p)),
\end{equation}
for $p\in F_{\beta g}$. The map $\psi \mapsto g \cdot \psi$
is semi-linear with respect to the action $S$,
$$ 
(g \cdot ({f \psi}))(p) = f(S(g^{-1})p) (g \cdot \psi)(p), 
$$
for $f\in C^\infty(F_{\alpha g})$. 
By allowing $g$ to range through the values of a bisection, we obtain
a semi-linear isomorphism of $\Ga E$: take any $b\in\mathcal{B}(G)$ 
and for $\psi\in\Ga E$ and $p\in E$ define
$$
(b\cdot \psi)(p) = (g\cdot \psi)(p)
$$
where $g = b((\beta\circ b)^{-1}(\phi(p))).$ It is clear that the resulting map
from $\mathcal{B}(G)$ to the group of semi-linear isomorphisms of $\Ga E$
is the same as that of Proposition \ref{prop:gpdslr}. This formulation perhaps
shows more clearly the effect of the action on sections of $E$, but it is
still unclear how to characterize the smoothness of $\rho$ in terms of the 
maps $\psi \mapsto g \cdot \psi$. 
\end{rmk}

\section{Appendix: some algebraic formulations}
As we remarked in the introduction, both the derivative endomorphisms and
the semi-linear isomorphisms are particular cases of the pseudo-linear 
transformations in the sense of Jacobson.
Here, we characterize the
pseudo-linear
endomorphisms of a left module ${\mathcal E}$
over a ring ${\mathcal C}$ in terms of twisted derivations.
In particular, the derivative endomorphisms 
are characterized in terms of derivations. We also characterize 
the semi-linear isomorphisms in terms of algebra automorphisms.
The case that is relevant to differential geometry is that of 
${\mathcal E}=\Gamma E$ 
and ${\mathcal C} = C^{\infty}(M)$, where 
$E$ is a vector bundle over a manifold $M$. For simplicity, we
shall formulate the results in this case only.

The notion of derivation of an algebra can be generalized as follows.
\begin{df}
Let $\mathcal A$ be an algebra and let $\alpha$ be a linear
endomorphism of $\mathcal A$.
A \emph{twisted derivation} of 
$\mathcal A$ with respect to $\alpha$
is a linear endomorphism $\mathcal U$ of $\mathcal A$ such that, 
for all $a, b \in {\mathcal A}$, 
\begin{equation} \label{twisted}
{\mathcal U}(ab) = {\mathcal U}(a)b+ \alpha(a){\mathcal U}(b). 
\end{equation}
\end{df}
It follows from \eqref{twisted} that, 
if ${\rm Im}({\mathcal U})$ has an element which is not a torsion element,
the map $\alpha$ is an algebra endomorphism.

\begin{df}
An $\mathbb R$--linear endomorphism $u$ of $\Gamma E$ is called 
\emph{pseudo-linear} if there exist 
${\mathbb R}$--linear endomorphisms $u^M$ and
$u_M$ of $C^{\infty}(M)$ such that
\begin{equation} \label{pseudo}
u(f \psi) = u^M(f) u(\psi) + u_M(f) \psi, 
\end{equation}
for all $f \in C^{\infty}(M)$ and $\psi \in \Gamma E$.
\end{df}

It follows from \eqref{pseudo} that, whenever 
there exists $\psi \in \Gamma E$ such that $u(\psi)$
and $\psi$ are linearly independent over the ring $C^{\infty}(M)$, 
(i) $u^M$ is an algebra endomorphism 
of $C^{\infty}(M)$, and (ii) $u_M$ is a twisted derivation of $C^{\infty}(M)$ 
with respect to $u^M$.

Clearly the derivative endomorphisms (resp., the semi-linear
isomorphisms)
of $\Gamma E$ are the
pseudo-linear  endomorphisms for which $u^M$ is the identity of 
$C^{\infty}(M)$ (resp., which are bijective and for which $u_M=0$).

Following Grothendieck \cite[16.5]{Grothendieck:1967},
we let ${\mathcal A}(E)$ be the vector space 
$C^{\infty}(M) \oplus \Gamma E$ 
with the commutative ${\mathbb R}$--algebra structure 
such that the product of any two elements of $C^{\infty}(M)$
is their product in the ring $C^{\infty}(M)$, the product of 
$f \in C^{\infty}(M)$ and $\psi \in \Gamma E$
is $f \psi$ in the 
$C^{\infty}(M)$--module $\Gamma E$, and the product of any two sections of
$\Gamma E$ vanishes. We can now formulate the following results,
whose proofs are all straightforward computations.

By a slight abuse of language, 
we let $u^M$ also denote the endomorphism of ${\mathcal A}(E)$
which coincides with $u^M$ on $C^{\infty}(M)$ and vanishes on 
$\Gamma E$.

\begin{prop}\label{prop3}
Let $u^M : C^{\infty}(M) \to C^{\infty}(M), ~ 
u_M : C^{\infty}(M) \to C^{\infty}(M)$ and $u : \Gamma E \to
\Gamma E$ be ${\mathbb R}$--linear maps.
Then  $u_M + u$ is a twisted derivation of the ${\mathbb R}$--algebra 
${\mathcal A}(E)$ with respect to $u^M$ if and only if
(i) $u_M$ is a twisted derivation of $C^{\infty}(M)$ with respect to $u^M$,
and (ii) $u$ is a pseudo-linear
endomorphism of $\Gamma E$ with associated endomorphisms of
$C^{\infty}(M)$, $u^M$ and $u_M$.
\end{prop}

When $u^M$ is the identity of $C^{\infty}(M)$, the twisted derivations
are ordinary derivations, so we obtain the following corollary.

\begin{prop}\label{prop1}
Let $u_M : C^{\infty}(M) \to C^{\infty}(M)$ and $u : \Gamma E \to
\Gamma E$ be ${\mathbb R}$--linear maps.
Then  $u_M+u$ is a derivation of the ${\mathbb R}$--algebra 
${\mathcal A}(E)$ if and only if
(i)~$u_M$ is a derivation of $C^{\infty}(M)$, and (ii) $u$ is a derivative
endomorphism of $\Gamma E$ with associated derivation of $
C^{\infty}(M)$, $u_M$.
\end{prop}

An analogous characterization holds for the semi-linear isomorphisms,
in terms of algebra automorphisms.

\begin{prop}\label{prop2}
Let $u^M : C^{\infty}(M) \to C^{\infty}(M)$ and $u : \Gamma E \to
\Gamma E$ be ${\mathbb R}$--linear isomorphisms.
Then  $u^M+u$ is an algebra 
automorphism of ${\mathcal A}(E)$ if and only if
(i) $u^M$ is an algebra automorphism of $C^{\infty}(M)$,
and (ii) $u$ is a semi-linear
isomorphism of $\Gamma E$ with associated automorphism of 
$C^{\infty}(M)$, $u^M$.
\end{prop}

The vector space of pseudo-linear endomorphisms of $\Gamma E$ is not closed 
under composition, but it follows from 
Proposition \ref{prop2} that the semi-linear isomorphisms 
of $\Gamma E$ form a group. In addition, just like the set of
twisted derivations, the set of pseudo-linear endomorphisms is not
closed under commutators, but it follows from Proposition \ref{prop1} that the
derivative endomorphisms of $\Gamma E$ form a Lie algebra.

\newcommand{\noopsort}[1]{} \newcommand{\singleletter}[1]{#1}

\end{document}